\documentclass[12pt]{article}
\usepackage[final]{epsfig}
\usepackage{graphics}
\usepackage{amsmath}
\usepackage{amsfonts}
\usepackage{latexsym}
\usepackage{amssymb}
\usepackage{graphicx}
\usepackage{url}

\newtheorem{theorem}{Theorem}
\newtheorem{lemma}[theorem]{Lemma}

\newcommand{\proofend}{$\Box$\bigskip}

\def\proof{\paragraph{Proof.}}

\begin{document}

\title{The Six Circles Theorem revisited}

\author{D. Ivanov\footnote{Moscow, Russia. e-mail: lesobrod@yandex.ru}
\ and 
S. Tabachnikov\footnote{
Department of Mathematics,
Penn State, 
University Park, PA 16802, 
and ICERM, Brown University, 
Box 1995, Providence, RI 02912.
e-mail: tabachni@math.psu.edu}
}

\date{}
\maketitle

\begin{figure}[hbtp]
\centering
\includegraphics[height=1.8in]{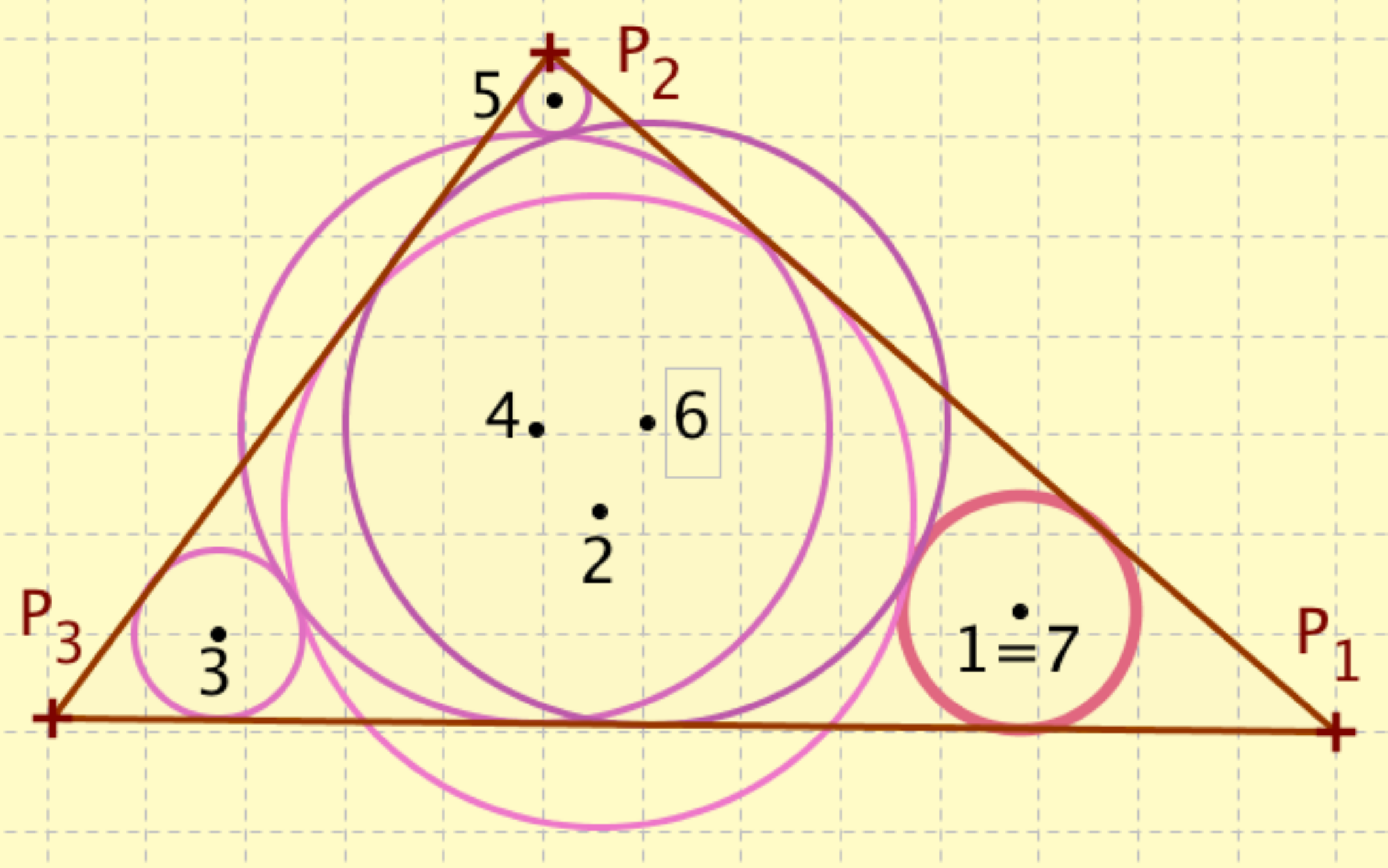}
\caption{The Six Circles Theorem: the centers of the consecutive circles are labeled $1,2,\dots,7$.}
\label{inside}
\end{figure}

{\bf Introduction}. Given a triangle $P_1 P_2 P_3$, construct a chain of circles: $C_1$, inscribed in the angle $P_1$; $C_2$, inscribed in the angle $P_2$ and tangent to $C_1$; $C_3$, inscribed in the angle $P_3$ and tangent to $C_2$; $C_4$, inscribed in the angle $P_1$ and tangent to $C_3$, and so on. The claim of The Six Circles Theorem is that this process is 6-periodic: $C_7=C_1$, see Figure \ref{inside}.

This beautiful theorem is one of many in the book \cite{EMT} which is a result of collaboration of three geometry enthusiasts, C. Evelyn, G. Money-Coutts, and J. Tyrrell. The following  is a quotation from John Tyrrell's obituary \cite{LS}:
\begin{quote}
 John also worked with two amateur mathematicians, C. J. A. Evelyn and G. B. Money-Coutts, who found theorems by using outsize drawing instruments to draw large figures. They then looked for concurrencies, collinearities, or other special features. The three men used to meet for tea at the Cafe Royal and talk about mathematics, and then go to the opera at Covent Garden, where Money-Coutts had a box.
\end{quote}
We refer to \cite{TP,Ri,FT,Ta,Tr} for various proofs and generalizations and to \cite{Ty} for a brief biography of C. J. A. Evelyn. See also \cite{Bo,wi,Wo} for Internet resources.
\medskip

{\bf A refinement}. The formulation of the Six Circles Theorem needs clarification. Firstly, there are two choices for each next circle; we assume that each time the smaller of the two circles tangent to the previous one is chosen (that is, the one which is closer to the respective vertex of the triangle). 
Secondly it well may happen that the next circle is tangent not to a side of the triangle but rather to its extension.

\begin{figure}[hbtp]
\centering
\includegraphics[width=5in]{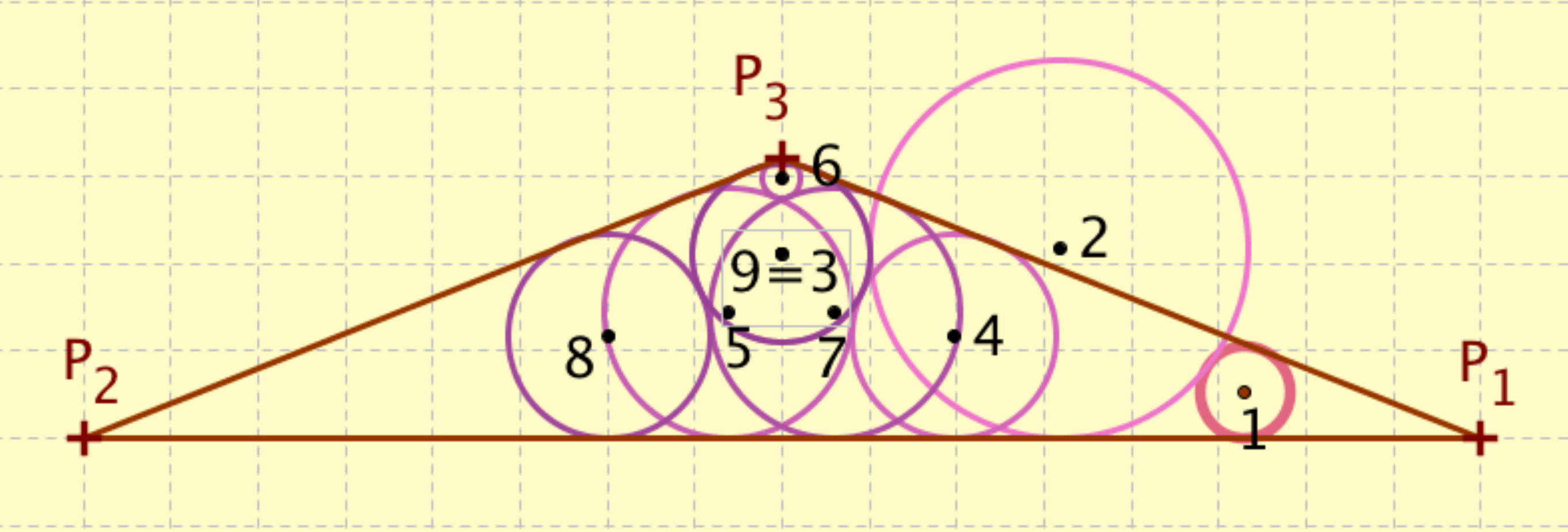}
\caption{The chain of circles is eventually 6-periodic with pre-period of length two: $C_{9}=C_3$, but $C_{8} \neq C_2$.}
\label{preper}
\end{figure}

The Six Circles Theorem, as stated at the beginning, holds for a chain of circles for which all tangency points lie on the sides of the triangle, not their extensions. And what about the latter case? Figure \ref{preper} shows what may happen.

\begin{theorem} \label{main}
Assume that, for the initial circle, at least one of the tangency points lies on a side of the triangle.
Then the chain of circles is eventually 6-periodic. One can choose the shape of a triangle and an initial circle so that the pre-period is arbitrarily long.
\end{theorem}

The existence of pre-periods is due to the fact that the map assigning the next circle to the previous one is not 1-1, that is, the inverse map is multi-valued. 

Concerning  the assumption that at least one of the tangency points of a circle with the sides of the angle of a tirangle
lies on a side of the triangle, and not its extension, we observe the following.

\begin{lemma} \label{tang}
If the first circle in the chain satisfies this assumption then so do all the consecutive circles. 
\end{lemma}

\begin{figure}[hbtp]
\centering
\includegraphics[width=3in]{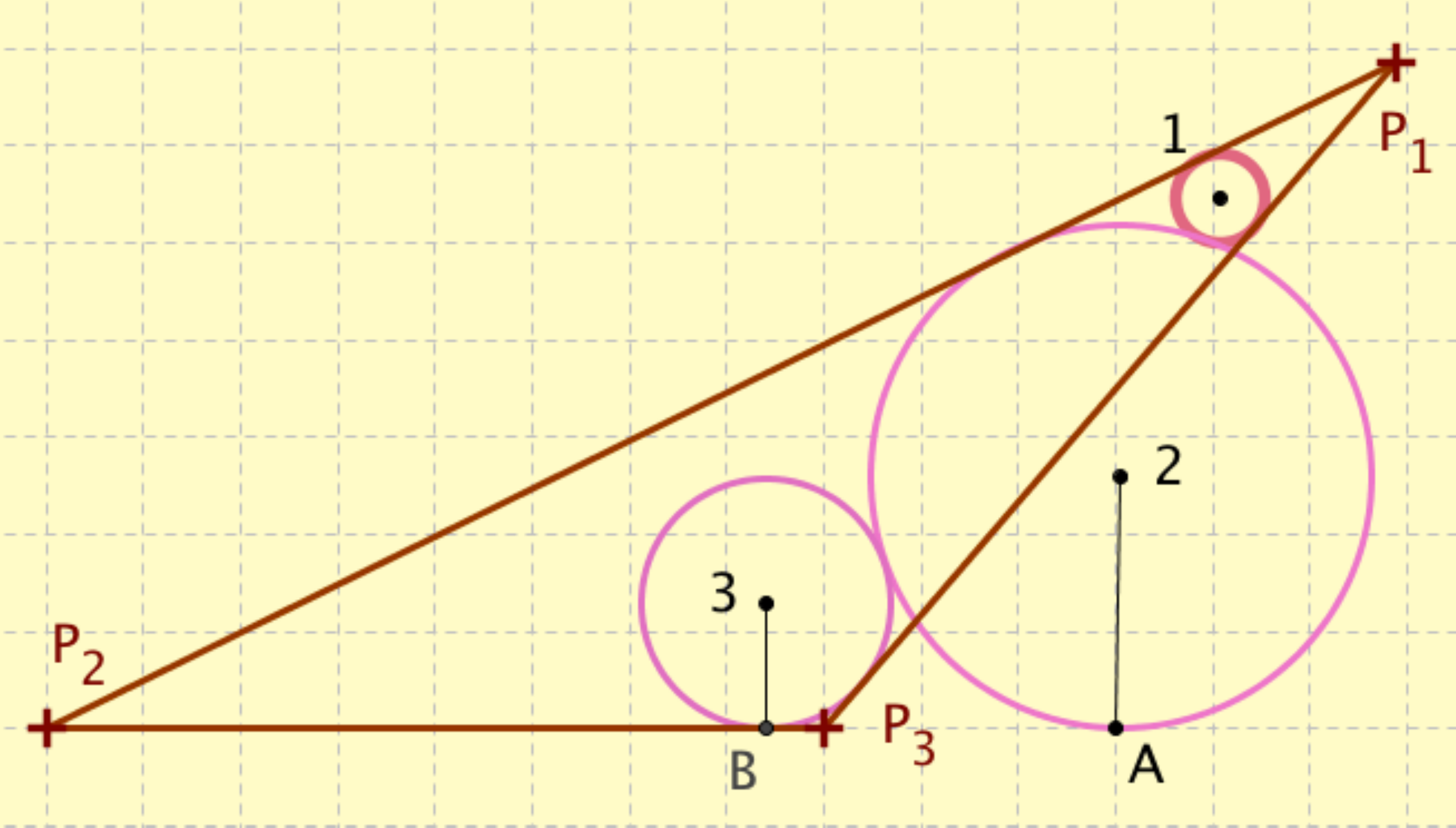}
\caption{The second case of equation (\ref{nextrad}): $|P_2 A| - |AB| + |BP_3| = |P_2 P_3|$.}
\label{case2}
\end{figure}

\proof
If circle $C_1$ touches side $P_1 P_2$ then circle $C_2$ also touches this side,  at a point closer to $P_2$ than the previous tangency point. Shifting the index by one, if circle $C_2$ does not touch side $P_2 P_3$ but 
touches side $P_1 P_2$ then it intersects side $P_1 P_3$, and the next circle $C_3$ touches side $P_1 P_3$, at a point closer to $P_3$ than the intersection points. See Figure \ref{case2} below for an illustration.
\proofend

\begin{figure}[hbtp]
\centering
\includegraphics[width=2.6in]{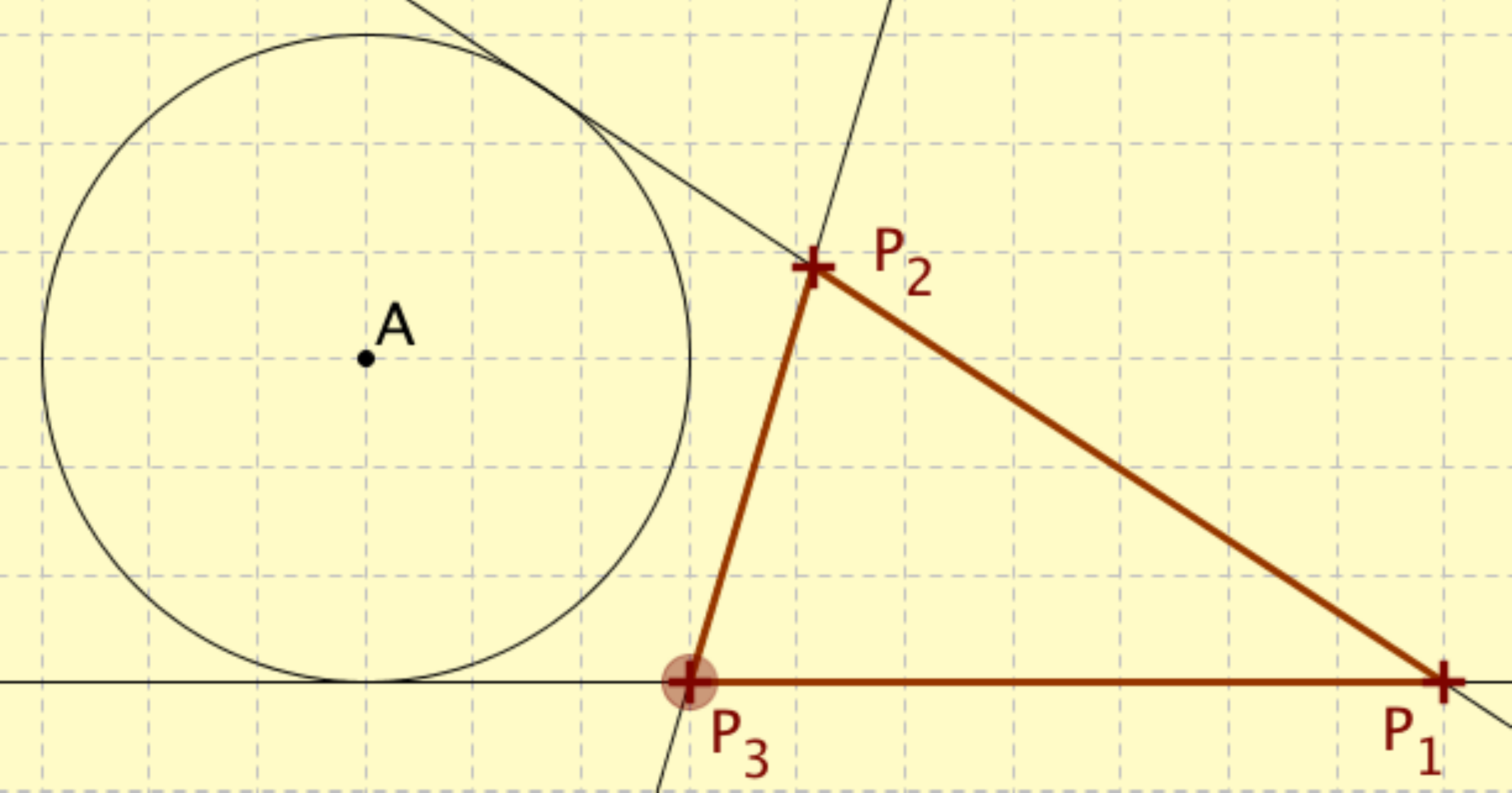} 
\includegraphics[width=2.6in]{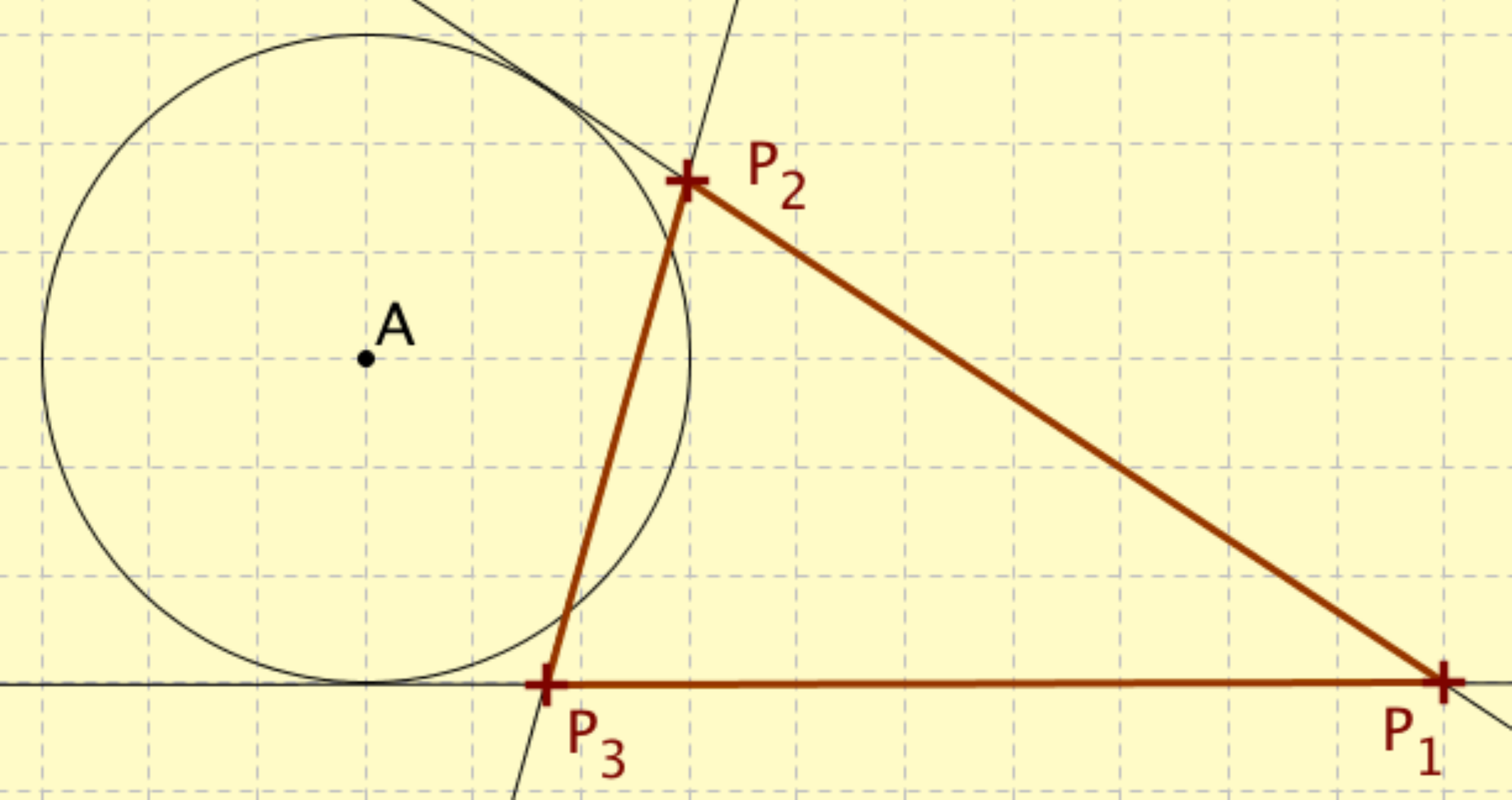}
\caption{When the initial circle touches the extensions of both sides of the triangle.}
\label{out}
\end{figure}

What about the case when the initial circle touches the extensions of both sides, $P_1 P_2$ and $P_1 P_3$? If the circle does not intersect  side $P_2 P_3$ then the next circle in the chain cannot be constructed, so this case is not relevant to us. If the first circle intersects  side $P_2 P_3$ then the next circle touches side $P_2 P_3$, and thus satisfies the assumption of Theorem \ref{main}, see Figure \ref{out}. Hence this assumption holds, starting with the second circle in the chain, and we may make it without loss of generality.

{\bf Beginning of the proof}. 
The proof consists of reducing the system to iteration of a piecewise linear function; this is achieved by a trigonometric change of variables (see \cite{EMT,Ta,Tr,FT} for versions of this approach). The choices of coordinates and various manipulations may look somewhat unmotivated; they are merely justified by the fact that they work.  The reader interested in a coordinate-free, but less elementary, approach is referred to \cite{Ta}.

Let us introduce notations. The angles of the triangle are $2\alpha_1, 2\alpha_2$ and $2\alpha_3$; its side  lengths are $a_1, a_2, a_3$ (with the usual convention that $i$th side is opposite $i$th vertex). Let $p=(a_1+a_2+a_3)/2$. We note that $p>a_i$ for $i=1,2,3$: this is the triangle inequality.
We denote the radii of the circles $C_i$ by $r_i,\ i=1,2,\ldots$ and assume that $C_i$ is  a circle that is inscribed into ($i$ mod $3$)-rd angle.

\begin{figure}[hbtp]
\centering
\includegraphics[width=3in]{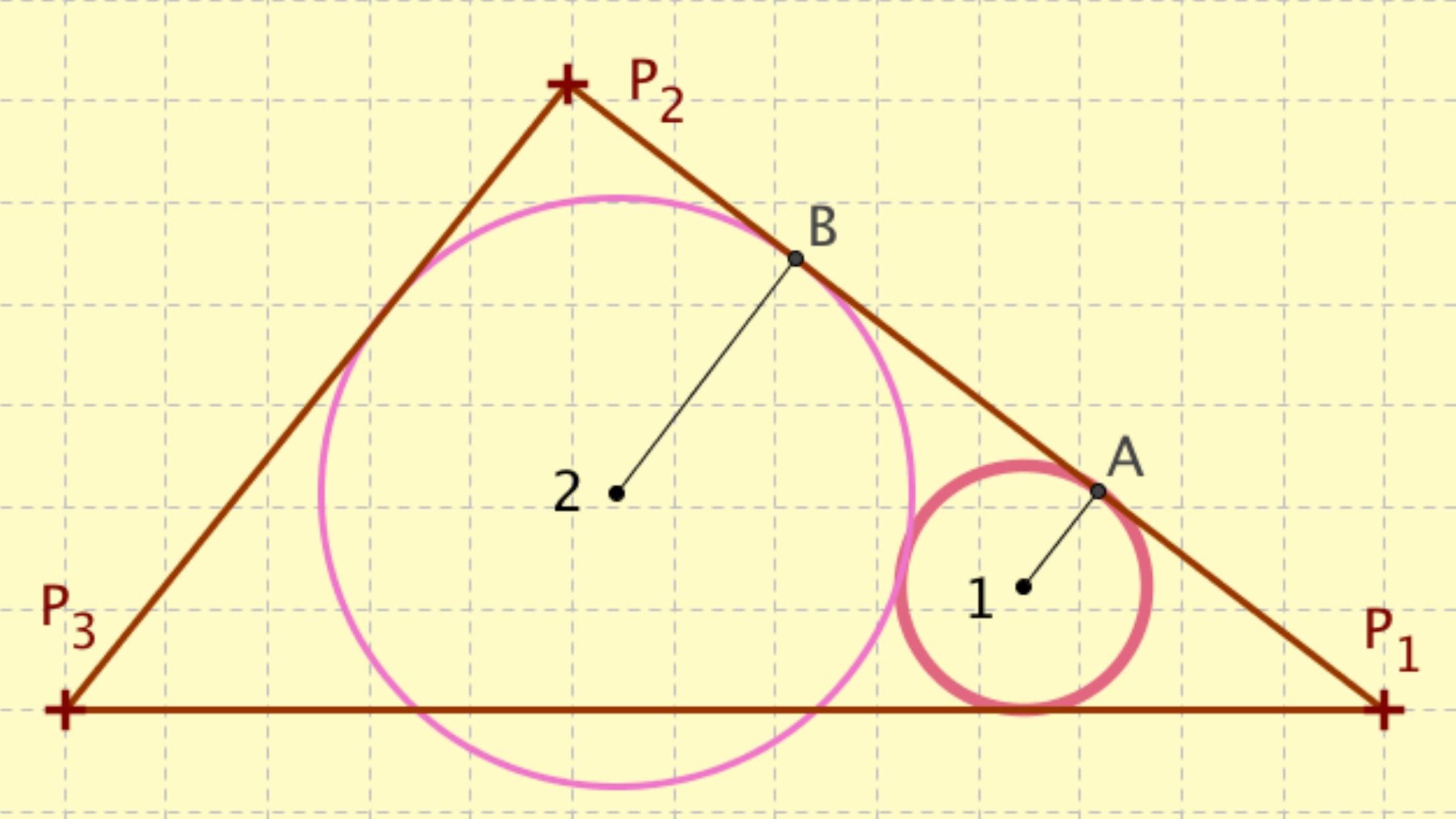} 
\caption{The first case of equation (\ref{nextrad}): $|P_1 A| + |AB| + |BP_2| = |P_1 P_2|$.}
\label{case1}
\end{figure}

\begin{figure}[hbtp]
\centering
\includegraphics[width=3in]{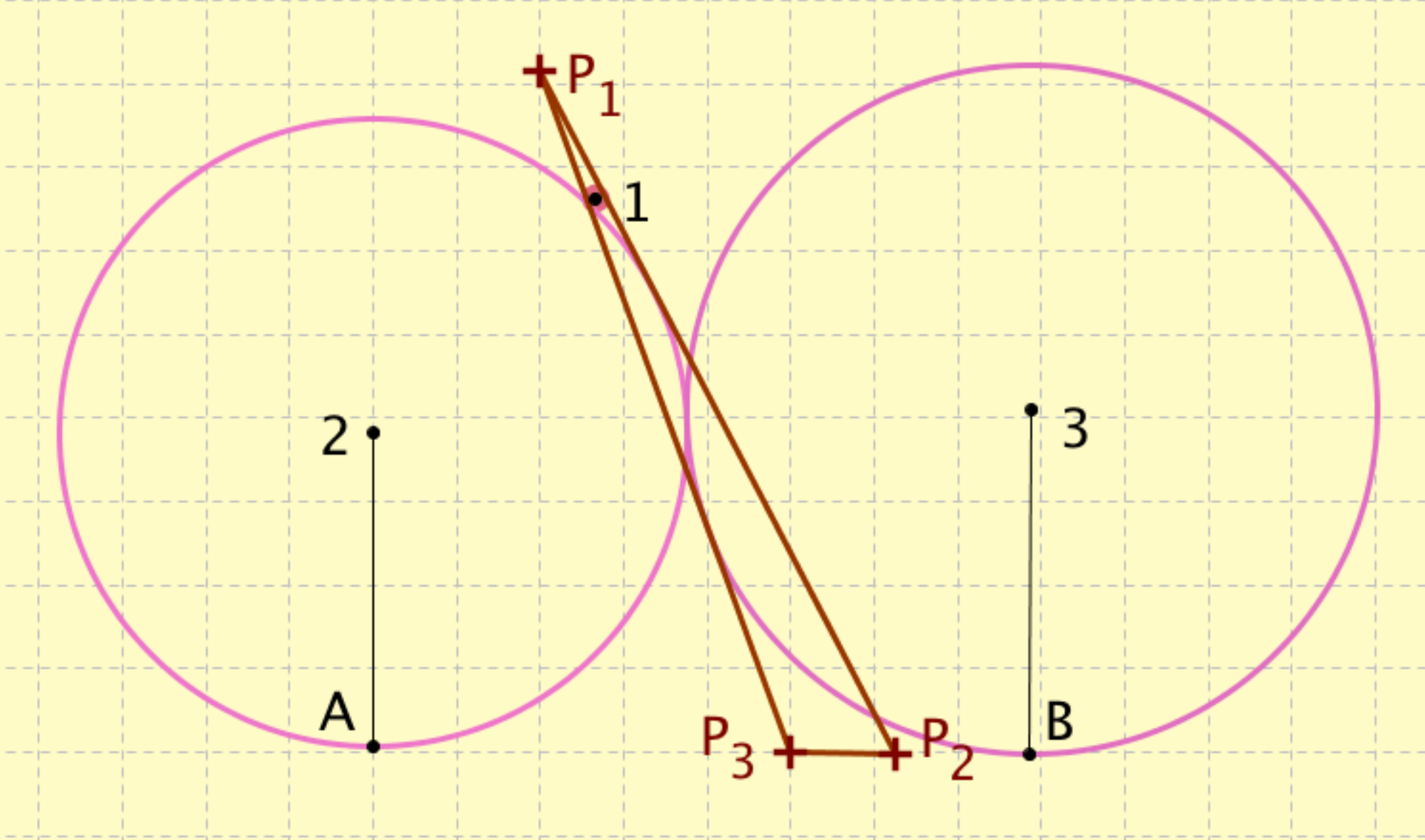}
\caption{Also the second case of equation (\ref{nextrad}): $|P_2 A| - |AB| + |BP_3| = |P_2 P_3|$.}
\label{minus}
\end{figure}

If two circles of radii $r_1$ and $r_2$ are tangent externally then the length of their common tangent segment (segment $AB$ in Figures \ref{case1}, \ref{case2}, \ref{minus}) is 
$$
\sqrt{(r_1+r_2)^2-(r_1-r_2)^2} = 2\sqrt{r_1 r_2}.
$$
Thus, depending on the mutual positions of the consecutive circles, as shown in Figures \ref{case1}, \ref{case2} and \ref{minus}, we obtain the equations
\begin{equation} \label{nextrad}
r_1\cot \alpha_1 + 2\sqrt{r_1 r_2} + r_2\cot \alpha_2=a_3\ \ {\rm or}\ \ r_1\cot \alpha_1 - 2\sqrt{r_1 r_2} + r_2\cot \alpha_2=a_3,
\end{equation}
or the cyclic permutation of the indices $1,2,3$ thereof. Specifically, if $C_1$ is tangent to the side $P_1P_2$ then we have the first equation (\ref{nextrad}), and if $C_1$ is tangent to the extension side $P_1P_2$ then we have the second equation. 
\medskip

{\bf Solving the equations}. 
Equations (\ref{nextrad}) determine the new radius $r_2$ as a function of the previous one, $r_1$. We shall solve these equations in two steps. First, introduce the notations
$$
u_1=\sqrt{r_1 \cot \alpha_1},\ \ e_3=\sqrt{\tan \alpha_1 \tan \alpha_2},
$$ 
and their cyclic permutations. Then (\ref{nextrad}) is rewritten as
\begin{equation} \label{inu}
u_1^2 \pm 2e_3u_1u_2+u_2^2=a_3,
\end{equation}
or 
\begin{equation} \label{expu}
u_1(u_1\pm e_3u_2) + u_2(u_2\pm e_3u_1)=a_3.
\end{equation}

Solve (\ref{inu}) for $u_2$:
\begin{equation} \label{solve}
u_2=-e_3u_1+\sqrt{a_3-(1-e_3^2)u_1^2},\ \ {\rm or}\ \ u_2=e_3u_1-\sqrt{a_3-(1-e_3^2)u_1^2},
\end{equation}
according as the sign in (\ref{inu}) is positive or negative. The minus sign in front of the radical in the second formula (\ref{solve}) is because our construction  chooses the smaller of the two circles tangent to the previous one. Likewise,  solve for $u_1$:
$$
u_1=-e_3u_2+\sqrt{a_3-(1-e_3^2)u_2^2},\ \ {\rm or}\ \ u_1=e_3u_2+\sqrt{a_3-(1-e_3^2)u_2^2},
$$
again depending on the sign in (\ref{inu}). The plus sign in front of the radical in the second formula is due to the fact that, going in the reverse direction, from $C_2$ to $C_1$, one chooses the greater of the two circles.
Substitute to (\ref{expu}) to obtain 
\begin{equation} \label{ufin}
u_1 \sqrt{a_3-(1-e_3^2)u_2^2} \pm u_2 \sqrt{a_3-(1-e_3^2)u_1^2} =a_3.
\end{equation}
The sign depends on whether $u_1^2$ is smaller or greater than $a_3$ (and if $u_1^2=a_3$ then $u_2=0$ in (\ref{solve})).
\medskip

{\bf Trigonometric substitution}. 
We shall rewrite the previous formula as the formula for sine of the sum or difference of two angles. To do so, we need a lemma. 

Given a triangle $ABC$, let $a,b,c$ be its sides, $p$ its semi-perimeter, and $2\alpha, 2\beta, 2\gamma$  its angles.

\begin{lemma} \label{tri}
One has
$$
1-\tan\alpha \tan\beta=\frac{c}{p}.
$$
\end{lemma}

\begin{figure}[hbtp]
\centering
\includegraphics[height=1.5in]{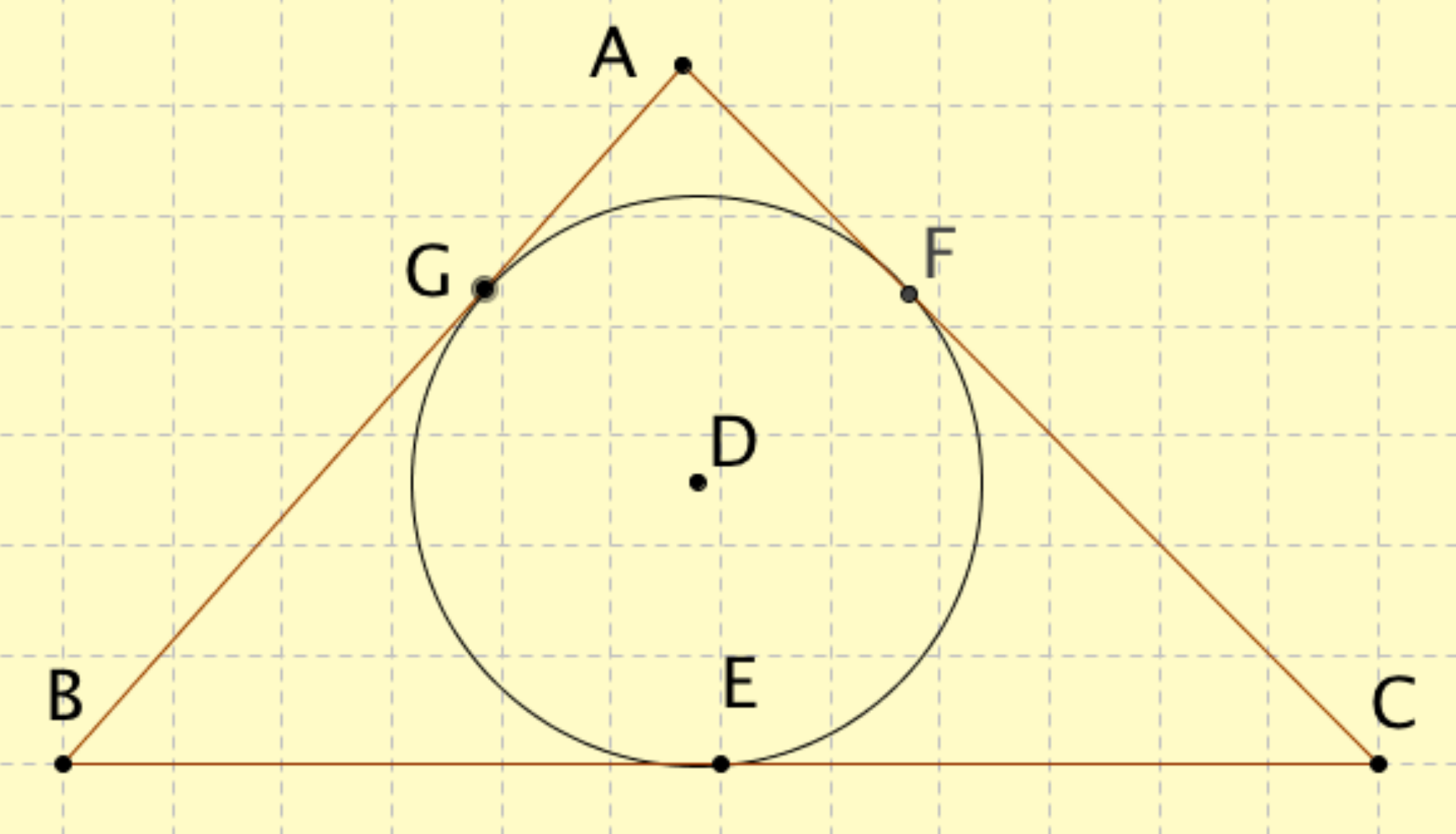}
\caption{To proof of Lemma \ref{tri}.}
\label{inscr}
\end{figure}

\proof
Let $R$ be the inradius and $S$ the area of the triangle. Let
$$
T_A=AF=AG,\ T_B=BG=BE,\ T_C=CE=CF,
$$ 
see Figure \ref{inscr}.
Then $p=T_A+T_B+T_C$, and $S=Rp$. By Heron's formula, $S=\sqrt{p T_AT_BT_C}$. Therefore
$R^2 p = T_AT_BT_C$.

On the other hand, $\tan \alpha=R/T_A, \tan\beta=R/T_B$, hence 
$$
1-\tan \alpha\tan\beta=1-\frac{R^2}{T_AT_B}=1-\frac{T_C}{p}=\frac{T_A+T_B}{p}=\frac{c}{p},
$$
as claimed.
\proofend

Using the lemma, we rewrite (\ref{ufin}) as
\begin{equation} \label{uff}
\frac{u_1}{\sqrt{p}}\sqrt{1-\frac{u_2^2}{p}}\pm \frac{u_2}{\sqrt{p}}\sqrt{1-\frac{u_1^2}{p}} = \sqrt{\frac{a_3}{p}}.
\end{equation}

We are ready for the final change of variables. Let
$$
\varphi_i = \arcsin\left(\frac{u_i}{\sqrt{p}}\right),\ \beta_i=\arcsin\left(\sqrt{\frac{a_i}{p}}\right).
$$

To justify the second formula, we note that $a_i < p$. Likewise, each circle is tangent to a side of the triangle, so $u_i^2$ is not greater than some side, and hence less than $p$. This justifies the first formula.

In the new variables, (\ref{uff}) rewrites as $\sin(\varphi_1 \pm \varphi_2) = \sin \beta_3$, where one has plus sign for $\varphi_1 < \beta_3$ and minus sign otherwise. Hence
\begin{equation} \label{module}
\varphi_2 = |\varphi_1-\beta_3|.
\end{equation}
This equation describes the dynamics of the chain of circles.

Before studying the dynamics of this function we note that the angles $\beta_i$ satisfy the triangle inequality, as the next lemma asserts. Assume that $\beta_1 \le \beta_2 \le \beta_3$.

\begin{lemma} \label{ineq}
One has $\beta_3 < \beta_1+\beta_2.$
\end{lemma}

\proof
We start by noting that $\sin \beta_i <1$ for  $i=1,2,3$, and that
$$
\sin^2 \beta_1+\sin^2 \beta_2+\sin^2 \beta_3=2\ \ {\rm or}\ \ \sin^2 \beta_3=\cos^2 \beta_1 + \cos^2 \beta_2.
$$
Assume that the triangle inequality is violated for some triangle. Since the inequality holds for an equilateral triangle, one can deform it to obtain a triangle for which $\beta_3 =\beta_1+\beta_2.$ Then
$$
\cos^2 \beta_1 + \cos^2 \beta_2=\sin^2 \beta_3 = (\sin \beta_1\cos\beta_2+\sin\beta_2\cos\beta_1)^2.
$$
It follows, after some manipulations, that 
$$
\sin \beta_1 \sin \beta_2 \cos\beta_1 \cos\beta_2=\cos^2\beta_1 \cos^2\beta_2\ \ {\rm or}\ \ \sin \beta_1 \sin \beta_2=\cos\beta_1 \cos\beta_2.
$$
Therefore 
$$
\cos(\beta_1+\beta_2)=0\ \ {\rm or}\ \ \beta_1+\beta_2=\frac{\pi}{2}.
$$
Hence $\sin^2 \beta_1+\sin^2 \beta_2=1$, and thus $\sin \beta_3 =1$. This is a contradiction.
\proofend
\medskip

{\bf Piecewise linear dynamics}. 
We are ready to investigate the function (\ref{module}). Although the dynamics of a piecewise linear function can be very complex \cite{Na}, ours is quite simple.

\begin{figure}[hbtp]
\centering
\includegraphics[height=1.4in]{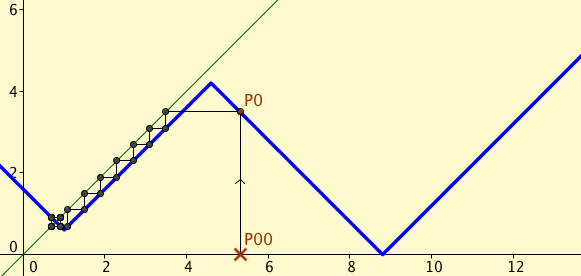}
\caption{Iteration of function $f(x)$ for $a=3.6, b=4.2$.}
\label{web}
\end{figure}

Iterating the map three times, with the values of the index $i=1,2,3$, yields the function
$y=|||x-\beta_1|-\beta_2|-\beta_3|$. We scale the $xy$ plane so that $\beta_1=1$ and rewrite the function as
\begin{equation} \label{funct}
f(x)=|||x-1|-a|-b|
\end{equation}
where $a\le b$ and $b<a+1$. We will show that every orbit of the map $f$ is eventually 2-periodic, see Figure \ref{web}.

\begin{figure}[hbtp]
\centering
\includegraphics[height=1.4in]{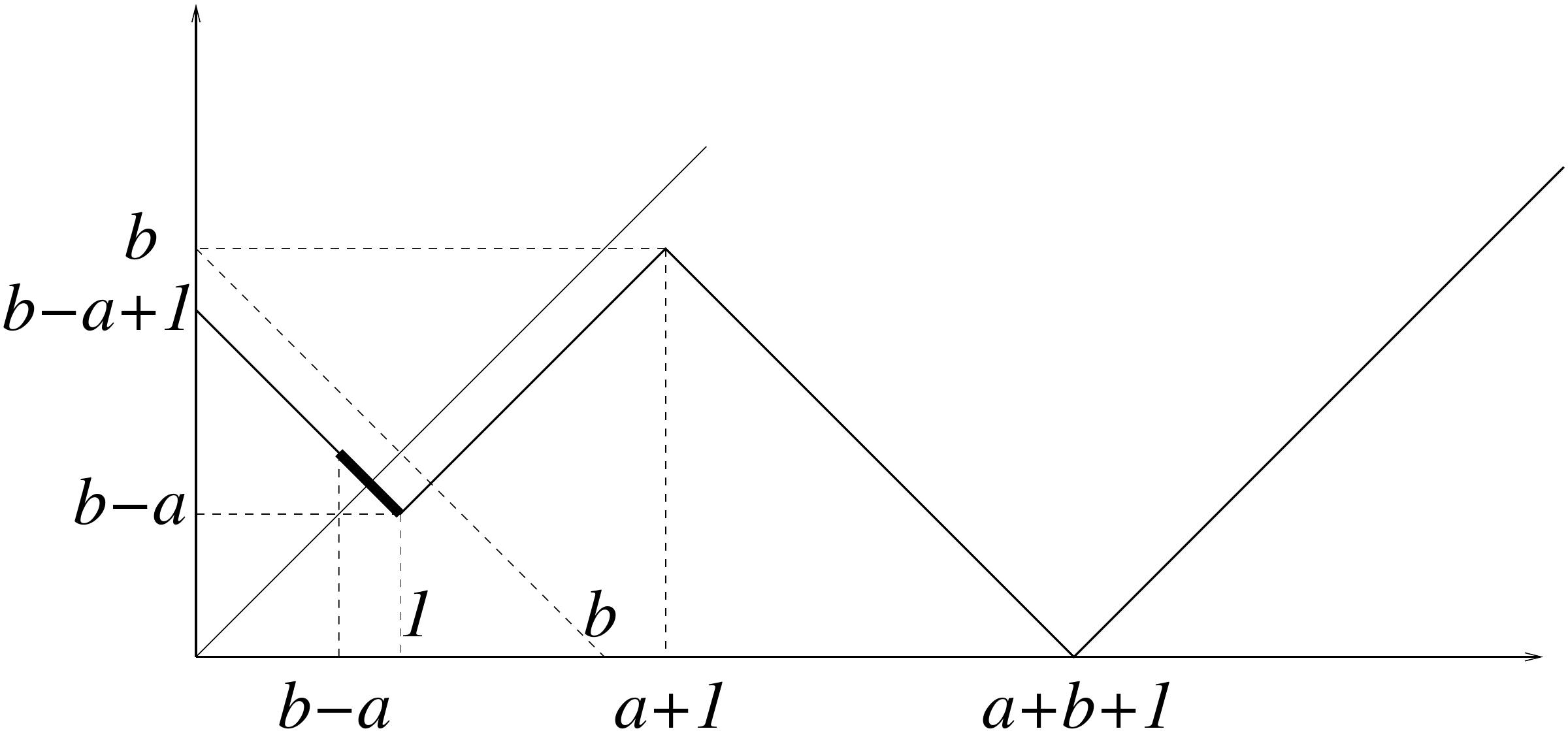}
\caption{The graph $y=f(x)$. The segment $[b-a,1]$ consists of 2-periodic points.}
\label{graph}
\end{figure}

The graph of $f(x)$ is shown in Figure \ref{graph} with  the characteristic points  marked.

It is clear that iterations of the function $f$ take every orbit to the segment $[0,b]$, and this segment is mapped to itself. Indeed, if $x\geq a+b+1$ then $f(x)=x-a-b-1$, and if $x \leq a+b+1$ then $f(x) \leq b$.
Thus iterations of the function $f$ will keep decreasing $x$ until it lands on $[0,b]$.

 Let 
$$
I_1=[0,b-a],\ I_2=[b-a,1],\ I_3=[1,b].
$$
Then $I_2$ consists of 2-periodic points, and we need to show that every orbit lands on this interval. Indeed, $f(I_1)=[1,b-a+1]\subset I_3$. On the other hand, each iteration of $f$ ``chops off" from the left a segment of length $1+a-b$ from $I_3$ and sends it 
to $I_2$. It follows that every orbit eventually reaches $I_2$. 

If $|I_2|=a+1-b$ is small, it may take an orbit a long time to reach $I_2$. For example, take $a=1$ and $b=2-\varepsilon$. Then, choosing $\varepsilon$ sufficiently small, one can make the pre-period of point $x=\varepsilon$ arbitrarily long. This choice corresponds to an isosceles triangle with the obtuse angle close to $\pi$ and a small initial circle $C_1$, compare with Figure \ref{preper}. 
\medskip

\begin{figure}[hbtp]
\centering
\includegraphics[width=3.5in]{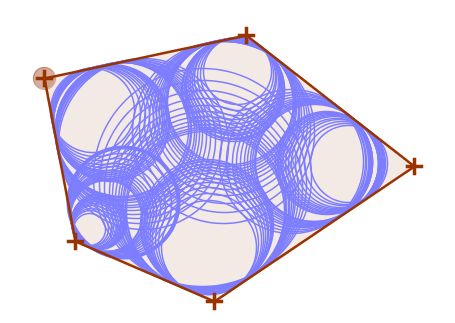}
\caption{A chain of circles in a pentagon.}
\label{penta}
\end{figure}

{\bf Final comments}. \\
 1) Although our considerations are close to those in \cite{EMT}, the authors of this book did not consider the pre-periodic behavior of the chain of circles. They addressed the issue of the two choices in each step of the construction and noted: 
\begin{quote} 
... we may make the first three sign choice quite arbitrarily provided that, thereafter, we make `correct' choices ...
\end{quote}
so that the chain becomes 6-periodic.\\
2) For a parallelogram, a similar phenomenon holds: the chain of circles is eventually 4-periodic but with a pre-period, see \cite{Tr}. Our analysis is similar to that of Troubetzkoy.\\
3) For $n>3$, the chain of circles inscribed in an $n$-gon is generically chaotic, see \cite{Tr} for a proof when $n=4$ and Figure \ref{penta} for an illustration when $n=5$. However, for every $n$, there is a class of $n$-gons enjoying $2n$-periodicity, see \cite{Ta}. Presumably, this periodicity is  also eventual, with an arbitrary long pre-period.\\

\begin{figure}[hbtp]
\centering
\includegraphics[width=4.5in]{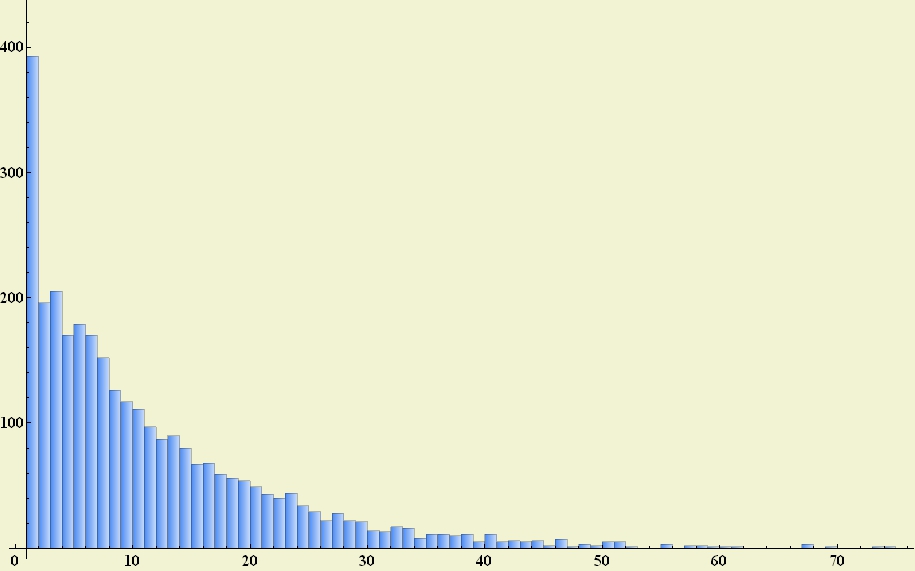}
\caption{The histogram represents 3000 chains of circles in a generic triangle. The selection, out of two, of each next circle in a chain is random. The horizontal axis represents the length of the pre-period, and the vertical the number of chains having this pre-period.}
\label{hist}
\end{figure}

4) A version of the Six Circles Theorem holds for curvilinear triangles made of arcs of circles \cite{EMT,TP,Ri}, and a generalization to $n$-gons is available as well \cite{Ta}. Again, one expects eventual periodicity with arbitrarily long pre-periods.\\
5) Constructing the chains of circles, we consistently chose the smaller  of the two circles tangent to the previous one. It is interesting to investigate what happens when other choices are made; for example,  one may toss a coin at each step. See Figure \ref{hist} for an experiment with a randomly chosen triangle.\\
6) The Six Circles Theorem is closely related with the Malfatti Problem: to inscribe three pairwise tangent circles into the three angles of a triangle; see, e.g., \cite{Gu} and the references therein. This 3-periodic chain of circles exists and is unique for every triangle; it corresponds to the fixed point of the function $f(x)$. See \cite{BZ} for a discussion of the Malfatti Problem close to our considerations.

\bigskip
{\bf Acknowledgments}. Most of the experiments that inspired this note and of the drawings were made in GeoGebra. The second author was supported by the NSF grant DMS-1105442. We are grateful to the referees for their criticism and advice.

\end{document}